\documentclass[10pt,a4paper]{article}
\usepackage[utf8]{inputenc}
\usepackage{amsmath}
\usepackage{amsfonts}
\usepackage{color}
\usepackage{amssymb}
\usepackage{tikz}
\usetikzlibrary{through}
\usetikzlibrary{intersections}
\usepackage{mathrsfs}
\usepackage{esint}
\usepackage[colorinlistoftodos]{todonotes}
\usepackage[colorlinks=true]{hyperref}
\hypersetup{
	colorlinks=blue,%
	citecolor=red,%
	filecolor=black,%
	linkcolor=blue,%
	urlcolor=blue
}
\usepackage[margin=1.3in]{geometry}

\usepackage[amsmath,thmmarks,hyperref]{ntheorem}
{
	\theoremstyle{nonumberplain}
	\theoremheaderfont{\bfseries}
	\theorembodyfont{\normalfont}
	\theoremsymbol{\mbox{$\Box$}}
	\newtheorem{pf}{Proof.}
}

\numberwithin{equation}{section}
\def\R{\mathbb{R}}
\def\S{\mathbb{S}}
\def\N{\mathbb{N}}

\def\e{\epsilon}
\newcommand{\ud}{\mathrm{d}}

\newtheorem{thm}{Theorem}[section]

\newtheorem{lem}{Lemma}[section]
\newtheorem{rem}{Remark}[section]

\newtheorem*{thm A}{Theorem A}
\newtheorem{cor}{Corollary}[section]
\usepackage{upgreek}
\usepackage{graphicx}
\newcommand{\DOI}[1]{\href{http://doi.org/#1}{#1}}
\allowdisplaybreaks
\title{article}
\begin{document}
	\title{\bf The sharp type Chern-Gauss-Bonnet integral and asymptotic behavior } 
\date{}
	\author{\medskip Shihong Zhang \\
		Nanjing University
	}
	
	\renewcommand{\thefootnote}{\fnsymbol{footnote}}
	\footnotetext[1]{S. Zhang: dg21210019@smail.nju.edu.cn}
	\maketitle

	%MS+++++++++++++++++++++ Abstract +++++++++++++++++++++++++
	{\noindent\small{\bf Abstract}: In this paper, we propose a sharp and quantitative criterion, which focuses solely on $Q$ curvature, to demonstrate the Chern-Gauss-Bonnet integral. In contrast to the previous results \cite{Chang Qing and Yang,Chang Qing and Yang invent,Hao}, we use a new approach that involves estimating the singular integral. Furthermore, we derive the asymptotic formula for the solution to the general $Q$ curvature equation.
 }

		\medskip 
		
		{{\bf $\mathbf{2020}$ MSC:} 58J05, 35J30}
		

		\section{Introduction}
		
		   The classic  Gauss-Bonnet formula states that suppose that $(M, g)$ is a compact manifold without a boundary, then
		\begin{align*}
			\int_{M}K_g\ud V_g=2\pi\chi(M),
		\end{align*}
		where $K_g$ represents the Gaussian curvature. In the case of complete surfaces, Huber \cite{Huber} and Finn \cite{Finn} demonstrated that if $K_g$ is absolutely integrable, then
		\begin{align}\label{G2}
			\int_{M}K_g\ud V_g\leq2\pi\chi(M)
		\end{align}
		and 
		\begin{align}\label{G3}
			\chi(M)-\frac{1}{2\pi}\int_{M}K_g\ud V_g=\sum_{i=1}^{k}\nu_i.
		\end{align}
		 Here $\nu_i$ is the isoperimetric ratio at the $i$-th end of $(M, g)$, i.e,
		 \begin{align*}
		 	\nu_i=\lim_{r_i\to\infty}\frac{|\partial B_{r_i}(0)|^2_g}{4\pi|B_{r_i}(0)|_g},
		 \end{align*}
		where $B_{r_i}(0)$ contains a compact set $ K_i$ and every $M\backslash K_i$ is diffeomorphic to $\R^2\backslash B_1(0)$.
		
		Later, Chang, Qing, and Yang \cite{Chang Qing and Yang,Chang Qing and Yang invent} made a successful extension of \eqref{G2} and \eqref{G3} to four dimensions. They showed that:
			\begin{thm A}[\protect{Chang, Qing and Yang \cite{Chang Qing and Yang}}]
			Suppose $(\R^4, e^{2u}|\ud x|^2)$ is complete with absolutely integrable $Q_ge^{4u}$, and supposed that 
			\begin{align}\label{scalar curvature assumption}
			R_g \mathrm{\,\,\, is\,\,\, nonnegative\,\,\, at\,\,\, infinity},
			\end{align}
			then
			\begin{align}\label{Q1}
			\int_{\R^4}Q_g\ud V_{g}\leq4\pi^2
			\end{align}
			and 
			\begin{align}\label{Q2}
			1-\frac{1}{4\pi^2}\int_{\R^4}Q_g\ud V_{g}=\lim_{r\to\infty}\frac{|\partial B_r(0)|_g^{\frac{4}{3}}}{4(2\pi)^{1/3}|B_r(0)|_g}.
			\end{align}
			\end{thm A}
		By studying the $n$-Laplace equation,	a general Huber type theorem have been established by S. Ma and Qing \cite{Ma and Qing}.

	The $Q$ curvature is a fascinating concept that has garnered significant attention in numerous papers. Some of the prominent references include \cite{Branson,Chang Gursky and Yang , Chang and Yang Ann, Fefferman,Fefferman&Hirachi}.
	
	In this context, we focus on the complete  manifold $(\R^n, g=e^{2u}|\ud x|^2)$ with  $Q_g\in L^1(\R^n,\ud V_g)$, where $n=2m\geq 4$, we also  denote
	\begin{align*}
	\alpha=\frac{1}{c_n}\int_{\R^n}Q_g\ud V_g=\frac{1}{c_n}\int_{\R^n}Q_g(x)e^{nu(x)}\ud x>0
	\end{align*}
	and
	\begin{align*}
	f(x)=Q_g(x)e^{nu(x)}\in L^1(\R^n).
	\end{align*}
		The $Q$ curvature in $\R^n$ is given by
	\begin{align}\label{key equation}
	\left(-\Delta\right)^m u=2Q_ge^{nu}	\qquad\qquad\mathrm{in}\quad\R^n.	
	\end{align}
	Furthermore, we use the notation  $o_{+}(1)$ to represent 
	\begin{align*}
	o_{+}(1)>0\qquad\mathrm{and}\qquad o_{+}(1)\to 0 \quad\mathrm{as}\quad |x|\to+\infty.
	\end{align*}		
	Let's define singular integral as $b(x)=\frac{1}{c_n}\int_{|x-y|<1}\log\frac{1}{|x-y|}f(y)\ud y$, the entire proof is focused on two estimations :
	\begin{enumerate}
		\item[(1)] The integral of $e^{b(x)}$ on a ray, namely,  $\int_{0}^{+\infty}e^{b(r,\theta_0)}\ud r<+\infty$. (Theorem \ref{main thm})
			\item[(2)] The upper bound of $b(x)$: for any $\e>0$, then $b(x)\leq \e\log|x|$ for $|x|\gg1$. (Theorem \ref{asy thm})
	\end{enumerate}

In this paper, we notice that  the result  \eqref{Q1} and \eqref{Q2} only involve $Q$ curvature. Thus, in my opinion,  it seems more appropriate to introduce additional assumptions that are specifically related to $Q$ curvature from a differential geometry perspective.
 We aim to generalize the results of Chang, Qing and Yang \cite{Chang Qing and Yang}  to high dimensions without scalar curvature assumption \eqref{scalar curvature assumption}. In the following, we will discuss some plausible conditions that only pertain to $Q_g$.
 
	\begin{enumerate}
		\item[(1)] $e^{-o_{+}(1)|x|^2}\leq Q_g(x)$ at infinity. The necessity of this assumption can be seen in the following example, which was proposed by Chang, Qing and Yang \cite{Chang Qing and Yang}. Let $u(x)=\log\frac{2}{1+|x|^2}+|x|^2$, then $g=e^{2u}|\ud x|^2$ is complete, but
		\begin{align*}
		\frac{1}{c_n}\int_{\R^n}Q_g\ud V_g=\frac{1}{c_n}\int_{\R^n}\left(-\Delta\right)^m u(x)\ud x=2>1.
		\end{align*}
		In the above example $Q_g(x)=e^{-n|x|^2}Q_{g_{\S^n}}$, which means that the decay of $Q_g$ does not seem to be too fast. In this paper, we prove that the decay rate
		$e^{-|x|^2}$ is sharp.  	In this paper, we assume that 
		\begin{align*}
		Q_g(x)\geq 0\qquad \mathrm{for}\qquad |x|>R_0.
		\end{align*}
		\item[(2)] $Q_ge^{nu}\in L^1(\R^n)$. This is a natural assumption for this problem. In fact, this condition will imply that $Q_g$ has some decay at infinity. Formally speaking, 
		\begin{align}\label{discuss u}
		u(x)\sim-\alpha \log|x| \qquad \mathrm{and}\qquad\alpha\leq1,
		\end{align}
		then \begin{align*}
		Q_g(x)e^{nu(x)}\sim\frac{Q_g(x)}{|x|^{n\alpha}}\in L^1(\R^n).
		\end{align*}
	
	\end{enumerate}

In \cite{Chang Qing and Yang}, their approach focused on analyzing the ODE in the radial symmetric case and subsequently reducing the general case to the radial case. In contrast, we adopt a completely different technique that involves proving the metric's normality (as seen in \eqref{normal metric condi}) through the completeness condition and the estimate of the singular integral.  Using this technique, we directly obtain \eqref{Thm 1.3.1}, then we also get \eqref{Thm 1.3.2} under the radial symmetric case and subsequently reduce the general case to the radial case. The core technique of the Theorem \ref{main thm} concentrate on  Lemma \ref{vital lem}, Lemma \ref{radial limit} and Lemma \ref{redc radial}. ``Estimating the length of a line near infinity'' is our primary tool for obtaining the following main theorem.
\begin{thm}\label{main thm}
	Let $(\R^n, g=e^{2u}|\ud x|^2)$ be a complete even dimensional manifold and  $Q_g\in L^1(\R^n,\ud V_{g})$. Suppose $e^{-o_{+}(1)|x|^2}\leq Q_g(x)$ at infinity, then
	\begin{align}\label{Thm 1.3.1}
	\int_{\R^n}Q_g\ud V_{g}\leq c_n
	\end{align}
	and
	\begin{align}\label{Thm 1.3.2}
	1-\frac{1}{c_n}\int_{\R^n}Q_g\ud V_{g}=\lim_{r\to\infty}\frac{|\partial B_r(0)|_g^{\frac{n}{n-1}}}{n\omega_{n-1}^{\frac{1}{n-1}}
		| B_r(0)|_g
	}.
	\end{align}
\end{thm}
\begin{rem}
	We note that if $R_g$ is nonnegative at infinity, then the inequality \eqref{Thm 1.3.1} has already been proven by H. Fang \cite{Hao}. We will consider some typical examples which satisfy all our  assumptions. Let
	\begin{align*}
	u_{\alpha}(x)=\frac{\alpha}{2}\log\frac{2}{1+|x|^2}\qquad\qquad\mathrm{for}\quad\alpha\in(0,1],
	\end{align*}
	then $g_{\alpha}=e^{2u_{\alpha}}|\ud x|^2$ is complete. Clearly,
	\begin{align*}
	Q_{g_{\alpha}}=\frac{1}{2}e^{-nu_{\alpha}}\left(-\Delta\right)u_{\alpha}=\frac{\alpha}{2} Q_{g_{\S^n}}\left(\frac{2}{1+|x|^2}\right)^{n\left(1-\frac{\alpha}{2}\right)}
	\end{align*}
	and
	\begin{align*}
	\int_{\R^n}Q_{g_{\alpha}}(x)e^{nu_{\alpha}(x)}\ud x=\frac{\alpha}{2}\omega_{n}Q_{g_{\S^n}}\leq c_n:=\frac{1}{2}\omega_{n}Q_{g_{\S^n}}.
	\end{align*}
	After a patient calculation, $u_{\alpha}$ will also satisfy  \eqref{Thm 1.3.2}.
\end{rem}
\begin{rem}
Under the assumptions stated in Theorem \ref{main thm}, if $Q_g\geq 0$ also holds in $\mathbb{R}^n$, it can be concluded that the scalar curvature is positive (for additional information, please refer to Lemma \ref{normal 1}).
\end{rem}
Recently, Wang \cite{Wang} made significant progress in the field of isoperimetric inequalities in higher dimensions. By utilizing a Sobolev inequality with weights, she demonstrated that assuming $(\R^n, e^{2u}|\ud x|^2)$ is complete and normal, i.e.,
\begin{align}\label{normal metric condi}
u(x)=\frac{1}{c_n}\int_{\R^n}\log\frac{|y|}{|x-y|}Q_g(y)e^{nu(y)}\ud y+C,
\end{align}
where $c_n=2^{n-2}\left(\frac{n-2}{2}\right)!\pi^{n/2}=\frac{1}{2}\omega_{n}Q_{g_{\S^n}}$ and $\omega_{n}=|\S^{n}|$. Moreover if 
\begin{align*}
\alpha_1:=\int_{\R^n}Q^{+}_g(y)e^{nu(y)}\ud y<c_n
\end{align*}
and 
\begin{align*}
\alpha_2:=\int_{\R^n}Q^{-}_g(y)e^{nu(y)}\ud y<+\infty,
\end{align*}
then there exist $C=C(n,\alpha_1,\alpha_2) $ such that for any bounded smooth domain $\Omega$ in $\R^n$,
\begin{align*}
|\Omega|_g\leq C(n,\alpha_1,\alpha_2)|\partial \Omega|_g^{\frac{n}{n-1}}.
\end{align*}
So, combining Wang's result with Theorem \ref{thm A1}, we have the following Corollary.
\begin{cor}
	With the same assumption as Theorem \ref{main thm} and suppose $\int_{\R^n}Q^{+}_g\ud V_{g}<c_n$, the isoperimetric inequality holds, i.e,  there exists $C(n,\alpha)$ such that for any bounded smooth domain $\Omega$ in $\R^n$
	\begin{align*}
	|\Omega|_g\leq C(n,\alpha)|\partial \Omega|_g^{\frac{n}{n-1}}.
	\end{align*}
\end{cor}

From the perspective of PDE, the above theorem can be stated as follows: under certain curvature assumptions, the solution to the PDE \eqref{key equation} exhibits specific asymptotic behavior. Specifically, under appropriate geometric assumptions, if one can prove that 
\begin{align}\label{more optimal asy formula}
	u(x)=-\alpha\log|x|+C+O\left(\frac{1}{|x|}\right),
\end{align}
then  \eqref{Thm 1.3.1} and \eqref{Thm 1.3.2} are trivial. However, in general, this does not appear to hold. In the second part, we will prove that $u(x)=-\alpha\log|x|+o(\log|x|)$ under very weak assumptions. To begin, we provide some historic background.  

In two-dimensional case, Cheng and Lin \cite{Cheng and Lin} showed that if $u$ satisfies
\begin{align*}
-\Delta u(x)=K(x)e^{2u(x)}
\end{align*}
with absolutely integrable $K(x)e^{2u(x)}$ and 
\begin{align}\label{chen lin decay}
	C_0e^{-|x|^{\beta}}\leq K(x)\leq C_0|x|^{\gamma}
\end{align}
 for some $\beta\in (0,1)$ and $\gamma>0$,  or $$C_0|x|^{-\gamma}\leq -K(x)\leq C_0|x|^{\gamma}$$for $|x|\gg 1$, then
\begin{align}\label{asym 1}
u(x)=\frac{1}{2\pi}\int_{\R^2}\log\left(\frac{|y|}{|x-y|}K(y)e^{2u(y)}\right)\ud y+C
\end{align}
and
\begin{align}\label{asym 2}
u(x)=-\alpha\log|x|+o(\log|x|) \qquad\mathrm{as}\qquad |x|\to+\infty.
\end{align}
Moreover, Cheng and Lin \cite{Cheng and Lin} have pointed out that $\beta<1$ is necessary for \eqref{asym 1} and \eqref{asym 2} without more assumptions.

In four dimension, if $Q_g$ is a constant in \eqref{key equation}, Lin \cite{Lin} creatively gave a systematic approach to study the higher order critical exponent elliptic equation. In high dimension,  the solutions are classified by Chang and Yang \cite{Chang and Yang} and Wei and Xu \cite{Wei and Xu} under $u(x)=o(|x|^2)$ constraint. Afterwards, L. Martinazzi \cite{Martinazzi} remove the assumption $u(x)=o(|x|^2)$. For more general $Q$, you can also refer \cite{Chen and Xu,Ge and Xu}. 

Our next goal is to develop more in-depth singular integral techniques and incorporate  elliptic estimates to confirm the predictions outlined in equations \eqref{discuss u}. Specifically, if we assume that $Q_g(x) \leq C_0|x|^{\gamma}$ as $|x|$ approaches infinity and $Q_g(x)\geq 0$, we can derive a high-dimensional version with respect to \eqref{asym 1} and \eqref{asym 2}. However, our condition allow for a faster decay of $Q_g$  than what is permitted in Cheng and Lin's assumption \eqref{chen lin decay}, and this is due to our completeness condition. It should be noted that the method used in \cite{Cheng and Lin} is no longer applicable in high dimensions, and proving Theorem \ref{asy thm} requires significant effort.
	  \begin{thm}\label{asy thm}
		With the same assumption as Theorem \ref{main thm}, we have
		\begin{align}\label{normal formu}
			u(x)=\frac{1}{c_n}\int_{\R^n}\log\frac{|y|}{|x-y|}Q_g(y)e^{nu(y)}\ud y+C,
		\end{align}
		and if $Q_g\geq 0$ and $Q_g(x)\leq |x|^{\gamma}$ at infinity for some $\gamma>0$, then
			\begin{align}\label{asy formula}
			u(x)=-\alpha\log|x|+o(\log|x|)\qquad\mathrm{as}\qquad |x|\to+\infty.
			\end{align}

	\end{thm}
\begin{rem}
In general, \eqref{asy formula} cannot be improved to \eqref{more optimal asy formula}, except for more restrictions on $Q_g$. We also give an assumption that it will derive \eqref{more optimal asy formula}, i.e, 
	\begin{align*}
		Q_g\left(\frac{x}{|x|^2}\right)|x|^{-2n+n\alpha}\in C^{\sigma}(B_1(0))
	\end{align*}
	for some $ \sigma\in (0,1)$.
\end{rem}

The paper is organized as follows. In Section \ref{sec 4.0}, we introduce some notations and basic estimates. In subsection \ref{sec 4.1}, we provide a detailed proof of the normal metric \eqref{normal formu} in high dimensions by estimating the singular integral. In subsection \ref{sec 5}, we establish the radial version of Theorem \ref{main thm} and then reduce the general case to the radial case. In Section \ref{final sec}, we derive the important asymptotic formula \eqref{asy formula}.

		\section{Proof of normal metrics } \label{sec 4} 
		
		\subsection{Preliminaries}\label{sec 4.0}
		In this paper, we often omit the volume form when it does not cause any confusion, for example,
	\begin{align*}
	\int_{B}f(y)\ud y:=\int_{B}f\qquad\mathrm{and}\qquad\int_{\partial B}f(y)\ud \sigma(y):=\int_{\partial B}f.
	\end{align*}
			For even number $n=2m$,	we consider the integral equation
		\begin{align}\label{v formula}
		v(x)=\frac{1}{c_n}\int_{\R^n}\log\frac{|x-y|}{|y|}f(y)\ud y,
		\end{align}
		where $c_n=2^{2m-2}\left(m-1\right)!\pi^{m}$, $f\in L^1(\R^n)$ and $f(x)\geq 0$ in $|x|>R_0$. Denote
		\begin{align*}
		\alpha=\frac{1}{c_n}\int_{\R^n}f(y)\ud y>0,
		\end{align*}
		obviously we have 
		\begin{align*}
		\left(-\Delta \right)^{m}v(x)=-2f(x)\quad\quad\quad\quad\mathrm{in}\quad \R^n.
		\end{align*}
	Here, we have gathered some well-known results (as seen in \cite{Lin}) regarding the singular integral, for their details you can see Appendix \ref{sec appendix} . Throughout the proof of normal metric (Theorem \ref{thm A1}), we will only provide detailed proofs for essential facts.
		
		\begin{lem}\label{upper bound}
			Suppose $v$ satisfies \eqref{v formula}  with $f\in L^1(\R^n)$ and $f(x)\geq 0$ in $|x|>R_0$ for some sufficiently big $R_0$ ,  if $|x|\gg1$ there holds
			\begin{align*}
			v(x)\leq \alpha\log|x|+C.
			\end{align*}
		
		\end{lem}

		\begin{lem}\label{lower bound}
			With the same assumption as Lemma \ref{upper bound}, then for any $\e>0$,  if $|x|\gg1$ we have
			\begin{align*}
			v(x)\geq(\alpha-\e)\log|x|-\frac{1}{c_n}\int_{|y-x|<1}\log\frac{1}{|x-y|}f(y)\ud y.
			\end{align*}
			
		\end{lem}

		\subsection{Normal metrics in even dimension}\label{sec 4.1}
		 For complete manifold $(\R^{n},g=e^{2u}|\ud x|^2)$, we consider the equation 
		\begin{align}\label{main equation for n dim}
		(-\Delta)^{m}u(x)=2Q_g(x)e^{nu(x)}\quad\quad\quad \qquad\mathrm{in}\qquad\R^n.
		\end{align}
		The main assumption is given by
		\begin{align}\label{main assumption for n dim}
			 e^{-o_{+}(1)|x|^2}\leq Q_g(x) \quad\mathrm{at \quad infinity}\qquad \mathrm{and}\qquad Q_g\in L^1(\R^n, \ud V_g).
		\end{align}

	We begin with mean value theorem and elliptic estimates for polyharmonic functions. The following two lemmas are well-known (see the appendix in \cite{Chen Lu Zhu}), so we also omit the details.
	\begin{lem}\label{Pizzetti}
	Suppose $\Delta^mh=0$ in $B_{2R}(x_0)$, then 
	\begin{align*}
		\fint_{B_{R}(x_0)}h(y)\ud y=\sum_{i=0}^{m-1}c_iR^{2i}\Delta^ih(x_0),
	\end{align*}
	where $c_i=\frac{n}{n+2i}\frac{(n-2)!!}{(2i)!!(2i+n-2)!!}$ for $i\geq1$ and $c_0=1$.
	\end{lem}

	\begin{lem}\label{high elliptic}
		If  $\Delta^m h=0$  in $B_4(0)$, then for any $\beta\in[0,1)$, $p\geq1$ and $k\geq1$, we have
		\begin{align*}
			||h||_{W^{k,p}(B_1(0))}\leq& C(k,p)||h||_{L^1(B_4(0))},\\
			||h||_{C^{k,\beta}(B_1(0))}\leq& C(k,p)||h||_{L^1(B_4(0))}.
		\end{align*}
	\end{lem}

	\begin{lem}\label{u+v=p}
	For $n\geq 4$, if $u$ solves \eqref{main equation for n dim} and satisfies \eqref{main assumption for n dim}, then 
	\begin{align}\label{formula a}
		u(x)+v(x)=p(x),
	\end{align}
	where $p(x)$ is a polynomial of degree at most $n-2$.
	\end{lem}	  
		\begin{pf}
			Let $h=u+v$, for any fixed $x_0\in\R^n$,  from the elliptic estimate of Lemma \ref{high elliptic}  we have
			\begin{align*}
			|\nabla^{n-1}h(x_0)|\leq& \frac{C}{R^{n-1}}\fint_{B_{R}(x_0)}|h|\\
			=&-\frac{C}{R^{n-1}}\fint_{B_{R}(x_0)}h	+\frac{2C}{R^{n-1}}\fint_{B_{R}(x_0)}h^{+}.
			\end{align*}
			Thanks to Lemma \ref{Pizzetti}, we know
			\begin{align*}
				\frac{C}{R^{n-1}}\fint_{B_{R}(x_0)}h=O(R^{-1}).
			\end{align*}
			Now we focus on the second part,  applying Lemma \ref{upper bound} and Jensen inequality we obtain
			\begin{align}\label{Lem 4.3.1}
				\frac{1}{R^{n-1}}\fint_{B_{R}(x_0)}h^{+}\leq& \frac{1}{R^{n-1}}\fint_{B_{R}(x_0)}u^{+}+C\frac{\log R}{R^{n-1}}\nonumber\\
				\leq& \frac{C}{R^{n-1}}\log\fint_{B_{R}(x_0)}e^{nu^{+}}+C\frac{\log R}{R^{n-1}}.
			\end{align}
			Since 
			\begin{align}\label{Lem 4.3.2}
				\int_{B_{R}(x_0)}e^{nu^{+}}\leq& \int_{B_{R}(x_0)}e^{nu}+1 \nonumber\\
				\leq&C(R_0,x_0) +CR^n+C\int_{B_{R}\backslash B_{R_0+|x_0|}(x_0)}e^{o_{+}(1)|x|^{2}}Q_ge^{nu}\nonumber\\
				\leq & C(R_0,x_0)+ CR^{n}+Ce^{CR^{2}}
			\end{align}
		for $R\gg1$, where $o_{+}(1)\to 0$ as $R_0\to+\infty$, but $o_{+}(1)\leq C$. Here $R_0$ is fixed, from \eqref{Lem 4.3.1} and \eqref{Lem 4.3.2} we arrive 
		\begin{align*}
			\frac{1}{R^{n-1}}\fint_{B_{R}(x_0)}h^{+}\leq \frac{C}{R^{n-3}}+\frac{C\log R}{R^{n-1}}\to 0
		\end{align*}
	 as $R\to\infty$.
		\end{pf}

	\begin{lem}\label{vital lem}
		Suppose  $(\R^n, g=e^{2u}|\ud x|^2)$ is complete and satisfies \eqref{main assumption for n dim}, then the following two statements holds.
		\begin{enumerate}
		\item[(1)] For $p(x)=p(r,\theta)$ in Lemma \ref{u+v=p}, there is no $\theta\in\S^{n-1}$ such that 
		\begin{align}\label{vital p1}
		\lim_{r\to\infty}\frac{p(r,\theta)}{r^{k}}=C(k,\theta)>0\qquad\qquad\mathrm{for}\quad k\geq 2.
		\end{align}
		\item[(2)]For $p(x)=p(r,\theta)$ in Lemma \ref{u+v=p}, there is no $\theta\in\S^{n-1}$ such that 
		\begin{align}\label{vital p2}
		\lim_{r\to\infty}\frac{p(r,\theta)}{r^{s}}\leq -C(s,\theta)<0\qquad\qquad\mathrm{for}\quad s>0.
		\end{align}
		\end{enumerate}
	\end{lem}
\begin{pf}
	For $(1)$, if not. Thus, for $r:=|x|>R>R_0$ and $\frac{x}{|x|}=\theta_0\in \S^{n-1}$, we know 
	\begin{align*}
	p(y)\geq c_0r^{k}\quad\quad\quad \qquad\mathrm{for}\quad y\in B_{1/r^{n-3}}(x),
	\end{align*}
	this follows by $|\nabla p(y)|\leq C|y|^{k-1}\leq C|y|^{n-3}\sim C|x|^{n-3}$.
	Then, thanks to Lemma \ref{upper bound},  we know 
	\begin{align*}
	u(y)=&-v(y)+p(y)\geq c_0r^{k}-\alpha\log r-C\\
	\geq& \frac{c_0}{2}r^{k} \qquad\qquad\mathrm{for}\quad y\in B_{1/r^{n-3}}(x).
	\end{align*}
	But we get
	\begin{align*}
	\int_{\R^n\backslash B_{R}}Q_ge^{nu}\geq& \int_{R}^{+\infty}\int_{\partial B_r\cap B_{r^{3-n}}(x)}Q_ge^{nu}\ud\sigma \ud r\\
	\geq& C\int_{R}^{+\infty}\frac{e^{\frac{nc_0}{2}r^{k}-o_{+}(1)r^2
	}}{r^{(n-3)(n-1)}}\ud r=+\infty,
	\end{align*}  
	where $o_{+}(1)\to 0$ as $R\to+\infty$, we choose $R>R_0$ sufficiently big. 
This is a contradiction with $\int_{\R^n}Q_ge^{nu}<+\infty$.

For $(2)$, if it's wrong. There exist $\theta_0\in\S^{n-1}$ such that 
\begin{align}\label{Lem 4.4.10}
\lim_{r\to\infty}\frac{p(r,\theta_0)}{r^{s}}\leq -C(s,\theta_0)<0\qquad\qquad\mathrm{for}\quad s>0.
\end{align}
We consider the length of curve near infinity, denote
\begin{align*}
III=\int_{R+1}^{+\infty}e^{u(r,\theta_0)}\ud r,
\end{align*}
where $R\gg1$ and $R>R_0$. From Lemma \ref{lower bound}, we know
\begin{align}\label{Lem 4.4.11}
III\leq C\int_{R+1}^{+\infty}\frac{e^{p(r,\theta_0)}}{r^{(\alpha-\e)}}II(r)\ud r,
\end{align}
where 
\begin{align*}
II(r)=\exp\left(\frac{1}{c_n}\int_{|r\theta_0-y|<1}\log\frac{1}{|r\theta_0-y|}f(y)\ud y\right).
\end{align*}
Since for any $r>R+1$, then $y\in \R^n\backslash B_{R}$. We rewrite $II(r)$ term as
\begin{align*}
II(r)=\exp\left(\int_{\R^n\backslash B_{R}}\sigma(R)\chi_{|r\theta_0-y|<1}\log\frac{1}{|r\theta_0-y|}\frac{f(y)}{c
	_n\sigma(R)}\ud y\right),
\end{align*}
where $\sigma(R)=\frac{1}{c_n}||f||_{L^1(\R^n\backslash B_{R})}<\frac{1}{2}$.	For $d\nu(y)=\frac{f(y)}{||f||_{L^1(\R^n\backslash B_{R})}}\ud y$, we apply the Jensen inequality to get
	\begin{align}\label{Lem 4.4.21}
II(r)\leq& \int_{\R^n\backslash B_{R}}\exp \left(\sigma(R)\chi_{|r\theta_0-y|<1}\log\frac{1}{|r\theta_0-y|}\right)\frac{f(y)}{c_n\sigma(R)}\ud y.
\end{align}
Plug \eqref{Lem 4.4.21} into \eqref{Lem 4.4.11} and by the Fubini's Theorem, we conclude
\begin{align}\label{Lem 4.4.31}
III\leq C\int_{\R^n\backslash B_{R}}\left(III_1(y)+III_2(y)\right)\frac{f(y)}{||f||_{L^1(\R^n\backslash B_{R})}}\ud y,
\end{align}
where
	\begin{align*}
III_1(y)=&\int_{I_y\cap(R+1,+\infty)}\frac{e^{p(r,\theta_0)}}{r^{(\alpha-\e)}|r\theta_0-y|^{\sigma(R)}}\ud r,\\
III_2(y)=&\int_{(R+1,+\infty)\backslash I_y}\frac{e^{p(r,\theta_0)}}{r^{(\alpha-\e)}}\ud r,
\end{align*}
and
\begin{align*}
I_{y}=\{r\theta_0|R+1<r<+\infty\}\cap B_1(y).
\end{align*} 
For any fixed $y\in\R^n\backslash B_{R}$,
denoted the center of $I_{y}$ by $y^{*}$, i.e, $y^{*}\in\{r\theta_0|R+1<r<+\infty\}$ and $(y^{*}-y)\cdot \theta_0=0$.
Here is a picture that shows their location:
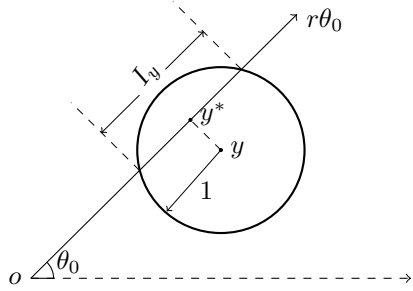
\begin{align*}
\begin{tikzpicture}
\draw[thick](2.5,1.7) circle (1.1);
\draw[fill] (2.5,1.7)  node [right] {$y$} circle(.02);% y coordiate
\coordinate[label=left: $o$] (o) at (0,0);
\node[below right] at (3.5,3.7) {$r\theta_0$};% r\theta_0
\draw[] (0,0)--(3mm,0pt)arc[start angle=0, end angle=45, radius=3mm] node[above=0.1pt, right]{$\theta_0$};%\theta_0
\draw[dashed, ->](0,0)--(5,0);
\draw[->](0,0)--(3.5,3.5);
\draw[very thin,dashed] (2.5,1.7) -- (2.1,2.1);
\draw[fill] (2.1,2.1)  node [right] {} circle(.02);
\node[above right] at (2.1,1.9) {$y^*$};
\draw[very thin,dashed] (1.43,1.43) -- (0.53,2.33);
\draw[very thin,dashed] (2.77,2.77) -- (1.87,3.67);
\draw[very thin,<-] (0.93,1.93) -- (1.43,2.43);
\draw[very thin,<-] (2.27,3.27) -- (1.77,2.77);
\node[above right] at (1.17,2.35) {\rotatebox{55}{$I_y$}};
\draw[very thin, ->] (2.5,1.7) -- (1.79,0.89);%radial 
\node[below right] at (2.1,1.4) {$1$};
\end{tikzpicture}
\end{align*}
Note that for any fixed $y\in\R^n$, $|r\theta_0-y^{*}|\leq|r\theta_0-y|$ and $|I_{{y}}|\leq2$, so we know
\begin{align}\label{Lem 4.4.41}
III_1\leq& \int_{(|y^{*}|-1, |y^{*}|+1)}\frac{e^{p(r,\theta_0)}}{r^{(\alpha-\e)}|r\theta_0-y^{*}|^{\sigma(R)}}\ud r\nonumber\\
\leq&C(\theta_0)\int_{(|y^{*}|-1, |y^{*}|+1)}\frac{1}{|r\theta_0-y^{*}|^{1/2}}\ud r
\leq C(\theta_0),
\end{align}
here we use $e^{p(r,\theta_0)}\sim e^{-cr^{s}}$ and $\sigma(R)\leq \frac{1}{2}$.
Clearly,
\begin{align}\label{Lem 4.4.51}
III_2=\int_{(R+1,+\infty)\backslash I_y}\frac{e^{p(r,\theta_0)}}{r^{(\alpha-\e)}}\ud r\leq C(\theta_0),
\end{align}
then \eqref{Lem 4.4.10}, \eqref{Lem 4.4.31}, \eqref{Lem 4.4.41} and \eqref{Lem 4.4.51} imply $III<+\infty$.
This is a contradiction with the completeness of $(\R^n, e^{2u}|\ud x|^2)$.

\end{pf}

	    \begin{thm}\label{thm A1}
	    		Suppose $(\R^n,g=e^{2u}|\ud x|^2)$ is complete
	    	with  $Q_g\in L^1(\R^n, \ud V_{g})$ . If $e^{-o_{+}(1)|x|^2}\leq Q_g(x)$ at infinity, then
	    	\begin{enumerate}
	    		
	    		\item[(1)]  the metric is normal, i.e,
	    		\begin{align}\label{nor condi}
	    		u(x)=\frac{1}{c_n}\int_{\R^n}\log\frac{|y|}{|x-y|}Q_g(y)e^{nu(y)}\ud y+C,
	    		\end{align}
	    		
	    		\item[(2)] 
	    		\begin{align}\label{ine a}
	    		\alpha:=\frac{1}{c_n}\int_{\R^n}Q_g(y)e^{nu(y)}\ud y\leq 1.
	    		\end{align}
	    	\end{enumerate} 
	    \end{thm}
    \begin{pf}
    	For $(1)$, let $p(x)=p_0+|x|p_1(\theta)+\cdots+|x|^kp_k(\theta)$, where $\theta=\frac{x}{|x|}\in \S^{n-1}$ and $k=\deg p\leq n-2$, we always assume $k\geq1$.
    	
    	 Step 1, we claim that $k$ is a even number. If $k$ is a odd number, there exists $\theta_0\in \S^{n-1}$   such that $p_k(\theta_0)<0$, thus
    	\begin{align*}
    	\lim_{r\to\infty}\frac{p(r,\theta_0)}{r^{k}}=p_k(\theta_0)<0.
    	\end{align*}
    	This cause a contradiction with \eqref{vital p2} in Lemma \ref{vital lem}. 
    	
     Step 2, $\sup_{\theta\in\S^{n-1}}p_k(\theta)\leq0$. Otherwise, there exists $\theta_0\in \S^{n-1}$   such that $p_k(\theta_0)>0$, hence
    	\begin{align*}
    	\lim_{r\to\infty}\frac{p(r,\theta_0)}{r^{k}}=p_k(\theta_0)>0,
    	\end{align*}
    	where $k\geq2$ due to Step 1. It's a contradiction with \eqref{vital p1} in Lemma \ref{vital lem}. 
    	
    	Step 3, $p\equiv C$. If not, clearly, $p_k(\theta)\not\equiv 0$, then there exists $\theta_0\in \S^{n-1}$   such that $p_k(\theta_0)<0$, 
    	\begin{align*}
    	\lim_{r\to\infty}\frac{p(r,\theta_0)}{r^{k}}=p_k(\theta_0)<0.
    	\end{align*}
    	This is also impossible by  \eqref{vital p2} in Lemma \ref{vital lem}. 
    	
    	For $(2)$, if $\alpha>1$, by the proof of Lemma \ref{vital lem}, we obtain 
    	 \begin{align*}
    	III\leq C\int_{\R^n\backslash B_{R}}\left(III_1(y)+III_2(y)\right)\frac{f(y)}{||f||_{L^1(\R^n\backslash B_{R})}}\ud y,
    	\end{align*}
    	where
    	\begin{align*}
    	III_1(y)=&\int_{I_y\cap(R+1,+\infty)}\frac{1}{r^{(\alpha-\e)}|r\theta_0-y|^{\sigma(R)}}\ud r,\\
    	III_2(y)=&\int_{(R+1,+\infty)\backslash I_y}\frac{1}{r^{(\alpha-\e)}}\ud r.
    	\end{align*}
    	We choose $R\gg1$, then $\e\ll 1$ , $\sigma(R)< 1/2$ and $\alpha-\e>1$. Note that,
    	\begin{align*}
    		\int_{-1}^{1}\frac{1}{|r|^{1/2}}\ud r <+\infty\qquad\mathrm{and}\qquad \int_{R+1}^{+\infty}\frac{1}{r^{\alpha-\e}}\ud r<+\infty.
    	\end{align*} 
    	We know $III<+\infty$, this is a contradiction with the completeness of $(\R^n, e^{2u}|\ud x|^2)$.
    	 \end{pf}

    We point out that \eqref{nor condi} and \eqref{ine a} will automatically imply the positivity of the scalar curvature under $Q_g\geq 0$ assumption.
      \begin{lem}\label{normal 1}
    	Suppose   the same condition as Theorem \ref{thm A1}  and $Q_g\geq 0$, then 
    	\begin{align*}
    	R_g(x)\geq-(n-1)\Delta u(x)e^{-2u(x)}>0.
    	\end{align*}
    \end{lem}
    \begin{pf}
    	By the conformal change of scalar curvature, we obtain
    	\begin{align}\label{Lem 4.6.1}
    	R_g(x)=&e^{-2u}(n-1)\left(-2\Delta u(x)-(n-2)|\nabla u(x)|^2\right)\nonumber\\
    	=&e^{-2u}(n-1)\left(2\Delta v(x)-(n-2)|\nabla v(x)|^2\right).
    	\end{align}
    	Clearly, we have
    	\begin{align}\label{Lem 4.6.2}
    	2&\Delta v(x)-(n-2)|\nabla v(x)|^2\nonumber\\
    	&=\frac{2(n-2)}{c_n}\int_{\R^n}\frac{f(y)}{|x-y|^2}\ud y-\frac{(n-2)}{c_n^2}\sum_{i=1}^{n}\left(\int_{\R^n}\frac{(x_i-y_i)}{|x-y|^2}f(y)\ud y\right)^2.
    	\end{align}
    	Notice that \eqref{ine a} implies that $\int_{\R^n}f(y)\ud y\leq c_n$, then 
    	\begin{align}\label{Lem 4.6.3}
    	\frac{(n-2)}{c_n}&\int_{\R^n}\frac{f(y)}{|x-y|^2}\ud y-\frac{(n-2)}{c_n^2}\sum_{i=1}^{n}\left(\int_{\R^n}\frac{(x_i-y_i)}{|x-y|^2}f(y)\ud y\right)^2\nonumber\\
    	&\geq\frac{(n-2)}{c_n^2}\int_{\R^n}f(y)\ud y\int_{\R^n}\frac{f(y)}{|x-y|^2}\ud y-\frac{(n-2)}{c_n^2}\sum_{i=1}^{n}\left(\int_{\R^n}\frac{(x_i-y_i)}{|x-y|^2}f(y)\ud y\right)^2\nonumber\\
    	&\geq 0.
    	\end{align}
    	The final step follows by Cauchy-Schwarz inequality. From \eqref{nor condi}, we arrive 
    	\begin{align}\label{Lem 4.6.4}
    	-\Delta u(x)=\frac{(n-2)}{c_n}\int_{\R^n}\frac{f(y)}{|x-y|^2}\ud y.
    	\end{align}
    	Combining with \eqref{Lem 4.6.1}, \eqref{Lem 4.6.2}, \eqref{Lem 4.6.3} and \eqref{Lem 4.6.4}, we finishing the proof.
    \end{pf}

		\subsection{Chern-Gauss-Bonnet formula in even dimension }\label{sec 5}
		In this section, we present a two-step proof of the Chern-Gauss-Bonnet formula. First, we assume that $u(x)=u(|x|)$, and we take a direct approach to prove it based on the normal condition stated in equation \eqref{integral f1} and the inequality given in  \eqref{inequality alpha 1}. Second, using Lemma \ref{asympotic formula}, we can reduce the general case to the radial case. Throughout this section, we make the following assumptions:
		 \begin{align}\label{integral f1}
		u(x)=\frac{1}{c_n}\int_{\R^n}\log\frac{|y|}{|x-y|}Q_g(y)e^{nu(y)}\ud y+C
		\end{align}
		and 
		\begin{align}\label{inequality alpha 1}
			\alpha=\frac{1}{c_n}\int_{\R^n}Q_g(y)e^{nu(y)}\ud y\leq 1.
		\end{align}
		
		\subsubsection{Radial symmetric case}
			\begin{lem}\label{radial limit}
			If $u(x)=u(|x|)$ and satisfies \eqref{integral f1}, then 
			\begin{align*}
			\lim_{r\to\infty}ru^{'}(r)=-\frac{1}{c_n}\int_{\R^n}f(y)\ud y.
			\end{align*}
		\end{lem}
		
		\begin{pf}
			Since 
			\begin{align*}
			ru^{'}(r)+\frac{1}{c_n}\int_{\R^n}f(y)\ud y=\frac{1}{c_n}\int_{\R^n}\frac{y\cdot(y-x)}{|y-x|^2}f(y)\ud y,
			\end{align*}
			we only need to proof $\int_{\R^n}\frac{y\cdot(y-x)}{|y-x|^2}f(y)\ud y\to 0$ as $|x|\to\infty$.
			By easy calculation, we know 
			\begin{align}\label{Re}
			I:=\int_{\R^n}\frac{y\cdot(y-x)}{|y-x|^2}f(y)\ud y=&\int_{0}^{+\infty}\int_{\partial B_s(0)}\frac{y\cdot(y-x)}{|y-x|^2}f(y)\ud\sigma(y)\ud r\nonumber\\
			=&\int_{0}^{+\infty}f(s)s\int_{\partial B_s(0)}\frac{y\cdot(y-x)}{|y||y-x|^2}\ud\sigma(y)\ud s\nonumber\\
			=&\int_{0}^{+\infty}f(s)s\int_{\partial B_s(0)}\frac{\partial \log|x-y|}{\partial \nu}\ud\sigma(y)\ud s.
			\end{align}
			We claim that for any  fixed $x$ and any $s>0$, we have 
			\begin{align}\label{Re 1}
			\int_{\partial B_s(0)}	\frac{\partial \log|x-y|}{\partial \nu}\ud\sigma(y)=\int_{ B_s(0)}	\Delta \log|x-y|\ud y=(n-2)\int_{ B_s(0)}\frac{1}{|x-y|^2}\ud y.
			\end{align}\begin{enumerate}
			
			\item[(1)] Firstly, if $x\not\in \bar{B}_s(0)$ the above identity is trivial by the Green's formula. 
			\item[(2)] Secondly, $x\in B_s(0)$, then for any $s-|x|>\rho>0$, we obtain
			\begin{align*}
			\int_{\partial B_s(0)}	\frac{\partial \log|x-y|}{\partial \nu}\ud\sigma(y)-\int_{\partial B_\rho(x)}	\frac{\partial \log|x-y|}{\partial \nu}\ud\sigma(y)
			=\int_{ B_s(0)\backslash B_\rho(x)}	\Delta \log|x-y|\ud y.
			\end{align*}
			Since 
			\begin{align*}
			\int_{\partial B_\rho(x)}	\frac{\partial \log|x-y|}{\partial \nu}\ud\sigma(y)=\omega_{n-1}\rho^{n-2}\to 0
			\end{align*}
		as $\rho\to 0$	and 
			\begin{align*}
			\int_{B_\rho(x)}	\Delta \log|x-y|\ud y=(n-2)\int_{ B_\rho(x)}\frac{1}{|x-y|^2}\ud y=\omega_{n-1}\rho^{n-2}\to 0.
			\end{align*}
			Hence, let $\rho\to 0$, we get the claim.
			\item[(3)] Finally, if $x\in \partial B_{s}(0)$, then there exists $\{x_n\}_{n=1}^{+\infty}$ such that $x_n\in B_s(0)$ and $\frac{x_n}{|x_n|}=\frac{x}{|x|}$. For any $\rho>0$, we denote
			\begin{align*}
				\partial B_{\rho}^{1}:=\partial B_s(0)\cap B_{\rho}(x)\qquad \mathrm{and}\qquad \partial B_{\rho}^2:=\partial B_{s}(0)\backslash \partial B_{\rho}^{1}.
			\end{align*}
			Then, we have 
			\begin{align}\label{Re 3}
				&\left|\int_{ \partial B_s(0)}\frac{\partial \log|x-y|}{\partial \nu}\ud\sigma(y)-\int_{ \partial B_s(0)}\frac{\partial \log|x_n-y|}{\partial \nu}\ud\sigma(y)\right|\nonumber\\
				&\leq C(n)\int_{ \partial B^1_\rho(x)}\frac{1}{|x-y|}+\frac{1}{|x_n-y|}\ud\sigma(y)+\int_{ \partial B_{\rho}^2(0)}\left|\frac{\partial \log|x-y|}{\partial \nu}-\frac{\partial \log|x_n-y|}{\partial \nu}\right|\ud\sigma(y)\nonumber\\
				&\leq C(n)\int_{\partial B^1_\rho(x)}\frac{1}{|x-y|}\ud\sigma(y)+o_{n}(1)=C(n)\rho^{n-2}+o_{n}(1),
			\end{align}
			where $o_n(1)\to 0$ as $n\to+\infty$ for any fixed $\rho$. And, 
			the last inequality follows from a fundamental result, i.e, $|x-y|\leq C|x_n-y|$. For the convenience of reader, we give its details. For any $y\in B_{\rho}(x)\cap \partial B_{s}(0)$, denote $\angle xoy=2\theta\ll \frac{\pi}{2}$, then $|x-y|=2|x|\sin\theta$.  We also  notice that $|x_n|\to |x|$, then
			\begin{align*}
				|x_n-y|=&\sqrt{|x|^2+|x_n|^2-2|x||x_n|\cos 2\theta}\\
				=&\sqrt{\left(|x|-|x_n|\right)^2+4|x||x_n|\sin^2\theta}\\
				\geq&\sqrt{2}|x|\sin\theta=\frac{\sqrt{2}}{2}|x-y|.
			\end{align*}
				
		\end{enumerate}
			And, 
			\begin{align}\label{Re 4}
				&\left|\int_{ B_s(0)}\frac{1}{|x-y|^2}-\frac{1}{|x_n-y|^2}\ud y\right|\nonumber\\
				&\leq\int_{ B_s(0)\backslash B_{\rho}(x)}\left|\frac{1}{|x-y|^2}-\frac{1}{|x_n-y|^2}\right|\ud y+\int_{ B_{\rho}(x)\cap B_{s}(0)}\left|\frac{1}{|x-y|^2}-\frac{1}{|x_n-y|^2}\right|\ud y\nonumber\\
				&\leq o_{n}(1)+\int_{B_{2\rho}(x)}\frac{1}{|x-y|^2}\ud y+\int_{B_{2\rho}(x_n)}\frac{1}{|x_n-y|^2}\ud y=o_{n}(1)+C\rho^{n-2}.
			\end{align} 
			From \eqref{Re 3} and \eqref{Re 4}, let $n\to+\infty$ firstly and $\rho\to 0$ secondly, then we complete the claim.
			
			Thus, combining \eqref{Re} with \eqref{Re 1}, we know 
			\begin{align*}
			I=(n-2)\int_{0}^{+\infty}f(s)s\int_{ B_s(0)}\frac{1}{|x-y|^2}\ud y\ud s.
			\end{align*}
			For any $\frac{1}{R_0}>\e>0$, let $|x|>\frac{1}{\e^2}$ we get
			\begin{align*}
			|I|=& \left|(n-2)\int_{0}^{\e|x|}+\int_{\e|x|}^{+\infty}f(s)s\int_{ B_s(0)}\frac{1}{|x-y|^2}\ud y\ud s\right|\\
			\leq&\frac{C}{(1-\e)^2|x|^2}\int_{0}^{\e|x|}|f(s)|s^{n+1}\ud s+C\int_{\e|x|}^{+\infty}f(s)s\int_{ B_s(0)}\frac{1}{|x-y|^2}\ud y\ud s\\
			\leq& \frac{C\e^2}{(1-\e)^2}\int_{0}^{+\infty}|f(s)|s^{n-1}\ud s+C\int_{1/\e}^{+\infty}f(s)s^{n-1}\ud s\\
			\leq& C\e^2||f||_{L^1(\R^n)}+C\int_{\R^n\backslash B_{1/\e}(0)}f(y)\ud y,
			\end{align*}
			where the third inequality follows by Rearrangement inequality  in Chapter 3 of \cite{Lieb and Loss},
			\begin{align*}
			\int_{ B_s(0)}\frac{1}{|x-y|^2}\ud y\leq	\int_{ B_s(0)}\frac{1}{|y|^2}\ud y\leq Cs^{n-2}.
			\end{align*}
			Let $\e\to 0$, we know $I\to 0$.
			
		\end{pf}

		\begin{lem}\label{radial gauss}
			With the same assumption as Lemma \ref{radial limit} and $(\R^n,e^{2u}|\ud x|^2)$ is complete, then
			 \begin{align}\label{Gauss Bonnet}
				1-\frac{1}{c_n}\int_{\R^n}Q_g(y)e^{nu(y)}\ud y=\lim_{r\to\infty}\frac{\left(\int_{\partial B_r}e^{(n-1)u(y)}\ud \sigma(y)\right)^{\frac{n}{n-1}}}{n\omega_{n-1}^{\frac{1}{n-1}}
					\left(\int_{B_{r}}e^{nu(y)}\ud y\right)
					}.
			\end{align}
				\end{lem}
			\begin{pf}
				For simplicity, we omit the volume form in the following proof, eg, $\int_{B_{r}}e^{nu(y)}\ud y:=\int_{B_{r}}e^{nu}$.  Applying Lemma \ref{upper bound} and \eqref{integral f1}, then for $|x|\gg1$ we have 
				 $$u(x)\geq-\alpha\log|x|-C.$$
				 Since $\alpha\leq1$, then
				\begin{align}\label{Lem 5.2.1}
					\int_{ B_r}e^{nu}\to\infty \qquad\mathrm{and}\quad \frac{\ud}{\ud r}	\int_{ B_r}e^{nu}=\int_{\partial B_r}e^{nu}>0.
				\end{align}
				Clearly,
				\begin{align*}
						&\frac{\frac{\ud}{\ud r}\left(\int_{\partial B_r}e^{(n-1)u}\right)^{\frac{n}{n-1}}}{n\omega_{n-1}^{\frac{1}{n-1}}
						\frac{\ud}{\ud r}\int_{B_{r}}e^{nu}
					}\nonumber\\
					=&\frac{ \frac{n}{n-1}\left(\int_{\partial B_r}e^{(n-1)u}\right)^{\frac{1}{n-1}}\left((n-1)\int_{\partial B_{r}}\frac{\partial u}{\partial r}e^{(n-1)u}+\frac{n-1}{r}\int_{\partial B_r}e^{(n-1)u}\right)  }{n\omega_{n-1}^{\frac{1}{n-1}}\int_{\partial B_{r}}e^{nu}}\nonumber\\
					=&\frac{e^{u(r)}\omega_{n-1}^{\frac{1}{n-1}}r\left(u^{'}(r)e^{(n-1)u(r)}+\frac{e^{(n-1)u(r)}}{r}\right)|\partial B_r|}{\omega_{n-1}^{\frac{1}{n-1}}|\partial B_r|e^{nu(r)}}
					=1+ru^{'}(r).
				\end{align*}
			 Now, from Lemma \ref{radial limit}, we get 
				\begin{align}\label{Lem 5.2.2}
						\lim_{r\to\infty}\frac{\frac{\ud}{\ud r}\left(\int_{\partial B_r}e^{(n-1)u}\right)^{\frac{n}{n-1}}}{n\omega_{n-1}^{\frac{1}{n-1}}
						\frac{\ud}{\ud r}\int_{B_{r}}e^{nu}
					}=1-\frac{1}{c_n}\int_{\R^n}Q_g(y)e^{nu(y)}\ud y.
				\end{align}
				Combining \eqref{Lem 5.2.1} with \eqref{Lem 5.2.2}, the L'Hopital's rule will imply that 
				\begin{align*}
						1-\frac{1}{c_n}\int_{\R^n}Q_g(y)e^{nu(y)}\ud y=\lim_{r\to\infty}\frac{\left(\int_{\partial B_r}e^{(n-1)u}\right)^{\frac{n}{n-1}}}{n\omega_{n-1}^{\frac{1}{n-1}}
						\left(\int_{B_{r}}e^{nu}\right)
					}.
				\end{align*}
			\end{pf}
	\begin{thm}\label{radial thm}
		Suppose that $u$ is radial function and  $Q_g\in L^1(\R^n, \ud V_{g})$.  If   $e^{-o_{+}(1)|x|^2}\leq Q_g(x)$ at infinity, then
		\begin{align*}
			\alpha=\frac{1}{c_n}\int_{\R^n}Q_g(y)e^{nu(y)}\ud y\leq 1
		\end{align*}
		and 
		\begin{align*}
			1-\frac{1}{c_n}\int_{\R^n}Q_g(y)e^{nu(y)}\ud y=\lim_{r\to\infty}\frac{|\partial B_r(0)|_g^{\frac{n}{n-1}}}{n\omega_{n-1}^{\frac{1}{n-1}}
				| B_r(0)|_g
			}.
		\end{align*}
	\end{thm}
	\begin{pf}
		See the Theorem \ref{thm A1} and Lemma \ref{radial gauss}.
	\end{pf}

	\subsubsection{General case}
	For $u\in C^{\infty}(\R^n)$, we define the spherical average
	\begin{align*}
		\bar{u}(r)=\fint_{\partial B_{r}(0)}u(y)\ud\sigma(y).
	\end{align*}
	The following Lemma \ref{asympotic formula} was originally proved by Chang, Qing, and Yang  \cite{Chang Qing and Yang} in the four-dimensional case. However, the result holds true for higher dimensions as well, so we omit the details of the proof.
		\begin{lem}[\protect{Chang Qing and Yang, \cite[Lemma 3.2]{Chang Qing and Yang}}]\label{asympotic formula}
		Suppose $\left(\R^n, e^{2u}|\ud x|^2\right)$ is normal metric. Then for any $k>0$,
		\begin{align*}
		\frac{1}{|\partial B_r|}\int_{\partial B_r}e^{ku}=e^{k\bar{u}(r)}e^{o(1)}.
		\end{align*}
	\end{lem}
\begin{pf}
	See the Page 523 of \cite{Chang Qing and Yang}.
\end{pf}

		\begin{lem}\label{redc radial}
			Suppose $(\R^n, g=e^{2u}|\ud x|^2)$ is complete  $Q_g\in L^1(\R^n, \ud V_{g})$ and $e^{-o_{+}(1)|x|^2}\leq Q_g(x)$ at infinity, then  $(\R^n, \bar{g}=e^{2\bar{u}(r)}|\ud x|^2)$ also satisfies the same properties. Moreover $\alpha_{\bar{g}}=\alpha_{g}$, namely, 
			\begin{align}\label{Lem 5.4.0}
				\frac{1}{c_n}\int_{\R^n}Q_{\bar{g}}(x)e^{n\bar{u}(x)}\ud x=\frac{1}{c_n}\int_{\R^n}Q_g(x)e^{nu(x)}\ud x.
			\end{align}
		\end{lem}
		\begin{pf}
			Since  $(\R^n, g=e^{2u}|\ud x|^2)$ is complete, then for any $\theta\in\S^{n-1}$,
			\begin{align*}
				\int_{0}^{+\infty}e^{u(r,\theta)}\ud r=+\infty.
			\end{align*}
			And we know $\fint_{\partial B_{r}(0)}e^{u(x)}\ud\sigma(x)=\fint_{\partial B_1(0)}e^{u(r,\theta)}\ud\sigma(\theta)$,
			by Fubini's theorem, then we know 
			\begin{align*}
				+\infty=\fint_{ \partial B_{1}(0)}\int_{0}^{+\infty}e^{u(r,\theta)}\ud r\ud\sigma(\theta)=\int_{0}^{+\infty}\fint_{ \partial B_{r}(0)}e^{u(x)}\ud\sigma(x)\ud r.
			\end{align*}
			From the  Lemma \ref{asympotic formula},  we get $\fint_{\partial B_{r}(0)}e^{u(x)}\ud\sigma(x)=e^{\bar{u}(r)}e^{o(1)}$, then 
			\begin{align}\label{Le 5.4.1}
				\int_{0}^{+\infty}e^{\bar{u}(r)}\ud r=+\infty,
			\end{align}
			i.e, $(\R^n, \bar{g}=e^{2\bar{u}(r)}|\ud x|^2)$ is complete. 
			Here it is apparent that
			\begin{align}\label{Le 5.4.2}
				2Q_{\bar{g}}(r)=&e^{-n\bar{u}(r)}\left(-\Delta\right)^m\bar{u}=e^{-n\bar{u}(r)}\fint_{ \partial B_{r}(0)}\left(-\Delta\right)^m u\nonumber\\
				=&e^{-n\bar{u}(r)}\fint_{ \partial B_{r}(0)}2Q_ge^{nu}.
			\end{align}
		Again through  Lemma \ref{asympotic formula}, we know 
			 $$\frac{1}{|\partial B_r|}\int_{\partial B_r}e^{nu}=e^{n\bar{u}(r)}e^{o(1)}.$$
			 	So, 
			 \begin{align}\label{Le 5.4.3}
			 Q_{\bar{g}}(r)\geq e^{-o_{+}(1)r^{2}}e^{-n\bar{u}(r)}\fint_{ \partial B_{r}(0)}e^{nu}\geq C_0e^{-o_{+}(1)r^{2}}\qquad\mathrm{for}\quad r\gg1.
			 \end{align}
			 For the last property, it follows by 
			 \begin{align}\label{Le 5.4.5}
			 	2\int_{\R^n}&Q_{\bar{g}}(x)e^{n\bar{u}(x)}\ud x=\int_{\R^n}\left(-\Delta\right)^m\bar{u}(x)\ud x\nonumber\\
			 	=&\int_{0}^{+\infty}\int_{\partial B_r}\left(-\Delta\right)^m\bar{u}\ud\sigma \ud r
			 	=\int_{0}^{+\infty}|\partial B_r| \left(-\Delta\right)^m\bar{u}\ud r\nonumber\\
			 	=&\int_{0}^{+\infty}|\partial B_r|\fint_{\partial B_{r}} \left(-\Delta\right)^mu\ud\sigma \ud r
			 	=\int_{\R^n}\left(-\Delta\right)^m u=2\int_{\R^n}Q_ge^{nu}.
			 \end{align}
			 	 With \eqref{Le 5.4.1}, \eqref{Le 5.4.2}, \eqref{Le 5.4.3} and \eqref{Le 5.4.5}, we complete the argument.
		\end{pf}
	\begin{rem}
		 It is worth noting that equation \eqref{Lem 5.4.0} is crucial in determining whether we can reduce the general case to the radial symmetric case. Specifically, the value of $\alpha$ in Theorem \ref{radial thm} and Theorem \ref{final thm} must be the same.
	\end{rem}

	\begin{thm}\label{final thm}
		Suppose $(\R^n, g=e^{2u}|\ud x|^2)$ is complete and  $Q_g\in L^1(\R^n, \ud V_{g})$ . Suppose $e^{-o_{+}(1)|x|^2}\leq Q_g(x)$ at infinity, then
		\begin{align*}
			1-\frac{1}{c_n}\int_{\R^n}Q_g(y)e^{nu(y)}\ud y=\lim_{r\to\infty}\frac{|\partial B_r(0)|_g^{\frac{n}{n-1}}}{n\omega_{n-1}^{\frac{1}{n-1}}
			| B_r(0)|_g
		}.
		\end{align*}
	\end{thm}
	\begin{pf}
		Similarly, by $\alpha\leq 1$ we know 
			\begin{align}\label{Thm 5.2.1}
		\int_{ B_r}e^{nu}\to\infty \qquad\mathrm{and}\quad \frac{\ud}{\ud r}	\int_{ B_r}e^{nu}=\int_{\partial B_r}e^{nu}>0.
		\end{align}
	Clearly,
		\begin{align}\label{Thm 5.2.2}
			&\frac{\frac{\ud}{\ud r}\left(\int_{\partial B_r}e^{(n-1)\bar{u}(r)}\right)^{\frac{n}{n-1}}}{n\omega_{n-1}^{\frac{1}{n-1}}
				\frac{\ud}{\ud r}\int_{B_{r}}e^{nu}
			}\nonumber\\
		=&	\frac{\left(\int_{\partial B_r}e^{(n-1)\bar{u}(r)}\right)^{\frac{1}{n-1}}\left(\int_{\partial B_{r}}\frac{\partial \bar{u}(r)}{\partial r}e^{(n-1)\bar{u}(r)}+\frac{1}{r}\int_{\partial B_r}e^{(n-1)\bar{u}(r)}\right)  }{\omega_{n-1}^{\frac{1}{n-1}}\int_{\partial B_{r}}e^{nu}}.
		\end{align}
		By the Lemma \ref{radial gauss}, we obtain
		\begin{align}\label{Thm 5.2.3}
			\lim_{r\to\infty}&
			\frac{\left(\int_{\partial B_r}e^{(n-1)\bar{u}(r)}\right)^{\frac{1}{n-1}}\left(\int_{\partial B_{r}}\frac{\partial \bar{u}(r)}{\partial r}e^{(n-1)\bar{u}(r)}+\frac{1}{r}\int_{\partial B_r}e^{(n-1)\bar{u}(r)}\right)  }{\omega_{n-1}^{\frac{1}{n-1}}\int_{\partial B_{r}}e^{n\bar{u}(r)}}\nonumber\\
			=&1-\frac{1}{c_n}\int_{\R^n}Q_g(y)e^{nu(y)}\ud y.
		\end{align}
			Lemma \ref{asympotic formula} implies that 
			\begin{align}\label{Thm 5.2.4}
				\lim_{r\to\infty}\frac{\int_{\partial B_{r}}e^{k\bar{u}(r)}}{\int_{\partial B_{r}}e^{ku}}=1 \qquad\mathrm{for}\quad k>0. 
			\end{align}
			Combining with \eqref{Thm 5.2.2}, \eqref{Thm 5.2.3} and \eqref{Thm 5.2.4}, we have
			\begin{align*}
				\lim_{r\to\infty}\frac{\frac{\ud}{\ud r}\left(\int_{\partial B_r}e^{(n-1)\bar{u}(r)}\right)^{\frac{n}{n-1}}}{n\omega_{n-1}^{\frac{1}{n-1}}
					\frac{\ud}{\ud r}\int_{B_{r}}e^{nu}
				}=1-\frac{1}{c_n}\int_{\R^n}Q_g(y)e^{nu(y)}\ud y.
			\end{align*}
			Hence,  \eqref{Thm 5.2.1} and L'Hopital's rule give that
			\begin{align}\label{Thm 5.2.5}
				\lim_{r\to\infty}\frac{\left(\int_{\partial B_r}e^{(n-1)\bar{u}(r)}\right)^{\frac{n}{n-1}}}{n\omega_{n-1}^{\frac{1}{n-1}}
					\int_{B_{r}}e^{nu}
				}=1-\frac{1}{c_n}\int_{\R^n}Q_g(y)e^{nu(y)}\ud y.
			\end{align}
		The final result follows by \eqref{Thm 5.2.4} and \eqref{Thm 5.2.5}.
		
	\end{pf}

Proof of Theorem \ref{main thm}: The follows by Theorem \ref{thm A1} and Theorem \ref{final thm}.

\section{Asymptotic behavior at infinity}\label{final sec}

Firstly, we recall the notation and basic facts:
\begin{align}\label{final sec v}
v(x)=\frac{1}{c_n}\int_{\R^n}\log\frac{|x-y|}{|y|}f(y)\ud y,
\end{align}
where $c_n=2^{2m-2}\left(m-1\right)!\pi^{m}$, and $f\in L^1(\R^n)$. Also, 
\begin{align*}
\alpha=\frac{1}{c_n}\int_{\R^n}f(y)\ud y,
\end{align*}
obviously we have 
\begin{align*}
\left(-\Delta \right)^{m}v(x)=-2f(x)\quad\quad\quad\quad\mathrm{in}\quad \R^n.
\end{align*}
In the following, we will assume that $f\geq 0$ !
Then, we also have the following facts:
\begin{align}\label{upper bounded}
v(x)\leq \alpha\log|x|+C,
\end{align}
and for any $\e>0$, if $|x|\gg1$
\begin{align}\label{lower bounded}
v(x)\geq(\alpha-\e)\log|x|-\frac{1}{c_n}\int_{|y-x|<1}\log\frac{1}{|x-y|}f(y)\ud y.
\end{align}
In the section, we focus on how to estimate the ``bad term'', i.e ,
\begin{align*}
	b(x)=\int_{|y-x|<1}\log\frac{1}{|x-y|}f(y)\ud y\leq \quad ?
\end{align*}
The whole process is somewhat difficult and highly nontrivial, in oder to do this , we need the following  basic estimates. The all efforts are concentrate on Lemma \ref{bad term}. During the process, we establish a very interesting  result (Lemma \ref{estimate of sphere integral}) involving singular integral. Specifically, for any $x\in \R^n$, there exists $C$ which is independent of $x$ such that 
	\begin{align*}
0<\fint_{\partial B_4(x)}-\left(-\Delta\right)^iv(y)\ud\sigma(y)\leq C,\qquad\mathrm{for}\quad i=1,\cdots,m-1.
\end{align*}
The basic structure is as follows: Lemma \ref{A 1 B} $\sim$ Lemma \ref{easy formula} will imply Lemma \ref{estimate of sphere integral};  Lemma \ref{estimate of sphere integral} $\sim$ Lemma \ref{improve exponential} finally cause Lemma \ref{bad term}.

\subsection{Estimate of singular integral}

\subsubsection{Weak harnack inequality of singular integral}

Now, we begin to consider the integral equation
\begin{align*}
v_k(x)=\int_{\R^n}\frac{f(y)}{|x-y|^k}\ud y.
\end{align*}
Firstly, we introduce the conception of $A_1$ weight. This is a nonnegative function $w\in A_1$ such that for any ball $B\subset\R^n$
\begin{align*}
\fint_{ B}w\leq C\inf_{ B}w
\end{align*}
or, equivalently, 
\begin{align*}
M(w)(x)=\sup_{r>0}\fint_{ B_{r}(x)}w\leq Cw(x)\qquad\qquad a.e.\quad x\in\R^n.
\end{align*}
For the basic knowledge of $A_p$ weight, you can refer the Chapter 5 of \cite{Stein}.

\begin{lem}\label{A 1 v}
	If $0<a<n$, then $\frac{1}{|x|^{a}}\in A_1$ weight.
\end{lem} 
\begin{pf}
	This is a easy excise, we omit the details.
\end{pf}

\begin{lem}\label{A 1 B}
	Suppose $f\geq 0$ and $f\in L^1(\R^n)$ with $0<k\leq n-2$, then $v_k(x)\in A_1$ weight and
	\begin{align}\label{A wieght 2}
	C\fint_{ B_{r}(x)}v_k\leq \inf_{ B_r(x)} v_k=\inf_{\partial B_r(x)} v_k\leq \fint_{\partial B_{r}(x)}v_k.
	\end{align}
\end{lem} 
\begin{pf}
	For any $r>0$, we have
	\begin{align*}
	\fint_{ B_{r}(x)}v_k(y)\ud y=&\fint_{ B_{r}(x)}\int_{\R^n}\frac{f(z)}{|y-z|^k}\ud z\ud y
	=\fint_{ B_{r}(x)}\frac{1}{|y-z|^k}\ud y\int_{\R^n}f(z)\ud z\\
	\leq&C\int_{\R^n}\frac{f(z)}{|x-z|^k}\ud z= Cv_k(x).
	\end{align*}
	The final inequality follows by  Lemma \ref{A 1 v}, i.e.,
	\begin{align*}
	\fint_{ B_{r}(x)}\frac{1}{|y-z|^k}\ud y=\fint_{ B_{r}(x-z)}\frac{1}{|y|^k}\ud y\leq \frac{C}{|x-z|^k}.
	\end{align*}
	Then, $M(v_k)(x)\leq Cv_k(x)$, i.e, $v_k\in A_1$. From easy calculation, we know 
	\begin{align*}
	-\Delta v_k(x)=k(n+k-2)\int_{\R^n}\frac{f(y)}{|x-y|^{k+2}}\ud y>0\qquad\mathrm{for}\quad 0<k<n-2.
	\end{align*} 
	And, 
	\begin{align*}
	-\Delta v_k=(n-2)\omega_{n-1}f\geq 0\qquad\mathrm{for}\quad k=n-2.
	\end{align*}
	Hence, \eqref{A wieght 2} follows by maximum principle and the definition of $A_1$ weight. 
\end{pf}

\subsubsection{More estimates of singular integral}

\begin{lem}\label{laplace of v}
	Suppose $v $ is given by \eqref{final sec v}, then 
	\begin{align*}
	\left(-\Delta\right)^k v(x)=\frac{d_k}{c_n}\int_{\R^n}\frac{f(y)}{|x-y|^{2k}}\ud y\qquad\quad\mathrm{for}\quad k=1,\cdots,m-1,
	\end{align*}
	where $d_{k+1}=2k(n-2k-2)d_k$ and $d_1=-(n-2)$.
\end{lem}
\begin{pf}
	The proof is very standard, we omit the details.
\end{pf}

In the following, we focus on estimating the term $\fint_{\partial B_4(x)}-\left(-\Delta\right)^iv(y)\ud \sigma(y)$, which will work in Lemma \ref{bad term}.
\begin{lem}\label{good integral}
	For $n>k>2$ and $x\in \R^n$, there holds 
	\begin{align*}
	\int_{\R^n\backslash B_4(x)}\frac{1}{|x-y|^{n-2}|y-z|^{k}}\ud y\leq C_1\chi_{B_2(x)}(z)+\frac{C_2\chi_{B^{c}_2(x)}(z)}{|x-z|^{k-2}}.
	\end{align*}
\end{lem}	
\begin{pf}
	For simplicity, we assume $x=0$ and let $I=\int_{\R^n\backslash B_4(0)}\frac{1}{|y|^{n-2}|y-z|^{k}}\ud y$.
	
	Case 1: $|z|\leq 2$.  Then $|y-z|>|y|-|z|\geq2$. Thus, we know 
	\begin{align*}
	|y|\leq |z|+|y-z|\leq 2|y-z|
	\end{align*}
	and 
	\begin{align*}
	|y-z|\leq|y|+|z|\leq \frac{3}{2}|y|.
	\end{align*}
	So, 
	\begin{align*}
	I\leq C\int_{\R^n\backslash B_4(0)}\frac{1}{|y|^{n+k-2}}\ud y\leq C_1.
	\end{align*}

	Case 2: $|z|\geq 2$, we split $\R^n$ into $\R^n=A_1\cup A_2 \cup A_3$, where
	\begin{align*}
	A_1=&\{y| |y-z|\leq|z|/2\},
	A_2=\{y| |y|\leq |z|/2\},\\
	A_3=&\{y||z|/2<|y|, |z|/2<|y-z|\}.
	\end{align*} We have
	\begin{align*}
	I\leq&\int_{A_1}+\int_{A_2}+\int_{A_3}\frac{1}{|y|^{n-2}|y-z|^{k}}\ud y\\
	:=&I_1+I_2+I_3.
	\end{align*}
	In $A_1$,  we know 
	\begin{align*}
	\frac{|z|}{2}\leq|y|\leq \frac{3}{2}|z|.
	\end{align*}
	So, 
	\begin{align*}
	I_1\leq \frac{C}{|z|^{n-2}}\int_{|y-z|<|z|/2}\frac{1}{|y-z|^k}\ud y\leq \frac{C}{|z|^{k-2}}.
	\end{align*} 
	Similarly, you can get 
	\begin{align*}
	I_2\leq \frac{C}{|z|^{k}}\int_{|y|<|z|/2}\frac{1}{|y|^{n-2}}\ud y\leq \frac{C}{|z|^{k-2}}.
	\end{align*}
	In $A_3$, you also know
	\begin{align*}
	\frac{1}{3}|y-z|\leq|y| \leq 3|y-z|.
	\end{align*}
	Thus, we obtain 
	\begin{align*}
	I_3\leq C\int_{|y-z|>|z|/2}\frac{1}{|y-z|^{n+k-2}}\ud y\leq C\frac{1}{|z|^{k-2}}.
	\end{align*}
\end{pf}

\begin{lem}\label{easy formula}
	For any $u\in L^1_{loc}(\R^n)$, then for any $R>0$ we have 
	\begin{align*}
	\int_{0}^{R}\frac{1}{|\partial B_r|}\int_{B_{r}}u(x)\ud x\ud r=\frac{1}{(n-2)\omega_{n-1}}\int_{B_{R}}\left(\frac{1}{|x|^{n-2}}-\frac{1}{R^{n-2}}\right)u(x)\ud x.
	\end{align*}
\end{lem}
\begin{pf}
	By element calculation, we know
	\begin{align*}
	\int_{0}^{R}\frac{1}{|\partial B_r|}\int_{B_{r}}u(x)\ud x\ud r=&\int_{0}^{R}\frac{1}{|\partial B_r|}\int_{0}^{r}\int_{\partial B_{s}}u(x)\ud\sigma \ud s\ud r\\
	=&\int_{0}^{R}\int_{\partial B_{s}}u(x)\ud\sigma \ud s\int_{s}^{R}\frac{1}{|\partial B_r|}\ud r\\
	=&\frac{1}{(n-2)\omega_{n-1}}\int_{B_{R}}\left(\frac{1}{|x|^{n-2}}-\frac{1}{R^{n-2}}\right)u(x)\ud x.
	\end{align*}
\end{pf}

\begin{lem}\label{estimate of sphere integral}
	If v is defined by \eqref{final sec v} and $f:=-\frac{1}{2}\left(-\Delta\right)^mv\in L^1(\R^n)$, then there exist $C$ which is independent of $x$ such that 
	\begin{align*}
	0<\fint_{\partial B_4(x)}-\left(-\Delta\right)^iv(y)\ud\sigma(y)\leq C,\qquad\mathrm{for}\quad i=1,\cdots,m-1,
	\end{align*}
	where $x\in\R^n$.
\end{lem}
\begin{pf}
	We argue it by induction. For $k=m-1$, we know
	\begin{align*}
	-\fint_{\partial B_{\rho}(x)}\frac{\partial}{\partial r}\left(-\Delta\right)^{m-1}v=\frac{1}{|\partial B_{\rho}(x)|}\int_{ B_{\rho}(x)}\left(-\Delta\right)^mv.
	\end{align*}
	Integral two sides from $0$ to $4$, then
	\begin{align*}
	LHS=&\left(-\Delta\right)^{m-1}v(x)-\fint_{\partial B_{4}(x)}\left(-\Delta\right)^{m-1}v\\
	=&\frac{1}{(n-2)\omega_{n-1}}\int_{\R^n}\frac{\left(-\Delta\right)^mv(y)}{|x-y|^{n-2}}\ud y-\fint_{\partial B_{4}(x)}\left(-\Delta\right)^{m-1}v,
	\end{align*}
	and by Lemma \ref{easy formula} we arrive
	\begin{align*}
	RHS=\frac{1}{(n-2)\omega_{n-1}}\int_{B_{4}(x)}\left(\frac{1}{|x-y|^{n-2}}-\frac{1}{4^{n-2}}\right)\left(-\Delta\right)^mv(y)\ud y.
	\end{align*}
	So, we get 
	\begin{align*}
	\fint_{\partial B_{4}(x)}&-\left(-\Delta\right)^{m-1}v\\
	=&\frac{1}{(n-2)\omega_{n-1}}\left(\int_{\R^n\backslash B_{4}(x)}\frac{-\left(-\Delta\right)^mv(y)}{|x-y|^{n-2}}\ud y+\frac{1}{4^{n-2}}\int_{B_{4}(x)}-\left(-\Delta\right)^mv(y)\ud y\right),
	\end{align*}
	we get the result for $k=m-1$. If we have know $k$ is right, we hope to prove $k-1$ is also true where $m-1\geq k\geq2$. Similarly, we know 
	\begin{align}\label{Lem 2.10.1}
	\fint_{\partial B_{4}(x)}&-\left(-\Delta\right)^{k-1}v\nonumber\\
	=&\frac{1}{(n-2)\omega_{n-1}}\left(\int_{\R^n\backslash B_{4}(x)}\frac{-\left(-\Delta\right)^kv(y)}{|x-y|^{n-2}}\ud y+\frac{1}{4^{n-2}}\int_{B_{4}(x)}-\left(-\Delta\right)^kv(y)\ud y\right).
	\end{align}
	By Lemma \ref{laplace of v} and Lemma \ref{good integral}, we conclude that
	\begin{align}\label{Lem 2.10.2}
	\int_{\R^n\backslash B_{4}(x)}\frac{-\left(-\Delta\right)^kv(y)}{|x-y|^{n-2}}\ud y=& C\int_{\R^n}\int_{\R^n\backslash B_{4}(x)}\frac{f(z)}{|x-y|^{n-2}|y-z|^{2k}}\ud y\ud z\nonumber\\
	\leq& C\int_{B_2(x)}f(z)\ud z+C\int_{\R^n\backslash B_2(x)}\frac{f(z)}{|x-z|^{2k-2}}\ud z
	\leq C.
	\end{align}
	For the second term, we know
	\begin{align*}
	-\left(-\Delta\right)^kv(y)=C\int_{\R^n}\frac{f(z)}{|y-z|^{2k}}\ud y=Cv_{2k}(y)\in A_1.
	\end{align*}
	By Lemma \ref{A 1 B}, we obtain 
	\begin{align}\label{Lem 2.10.3}
	\fint_{B_4(x)}-\left(-\Delta\right)^kv=&C\fint_{B_4(x)}v_{2k}
	\leq C\inf_{B_4(x)}v_{2k}=C\inf_{\partial B_4(x)}v_{2k}\nonumber
	\\
	\leq&C\fint_{\partial B_4(x)}v_{2k}=C\fint_{\partial B_4(x)}-\left(-\Delta\right)^{k}v\leq C,
	\end{align}
	where we use the induction for $k$. From \eqref{Lem 2.10.1} , \eqref{Lem 2.10.2} and \eqref{Lem 2.10.3}, we complete the argument.
\end{pf}

\begin{lem}\label{negative part L 1 estimate}
	If $f\in L^1(\R^n)$and $f\geq0$, then  for $R_0\gg1$ 
	\begin{align*}
	\int_{\R^n\backslash B_{R_0}(0)}v^{-}\leq C.
	\end{align*}
\end{lem}
\begin{pf}
	By lower bound estimate \eqref{lower bounded}, we know 
	\begin{align*}
	\int_{\R^n\backslash B_{R_0}(0)}v^{-}(x)\ud x\leq& C\int_{\R^n\backslash B_{R_0}(0)}\int_{|y-x|<1}\log \frac{1}{|x-y|}f(y)\ud y\ud x\\
	\leq&C\int_{\R^n}\int_{\R^n}\chi_{|x-y|<1}\log \frac{1}{|x-y|}f(y)\ud y\ud x\\
	\leq&C\int_{B_1(y)}\log \frac{1}{|x-y|}\ud x\int_{\R^n}f(y)\ud y\leq C.
	\end{align*}
\end{pf}

Finally, we would like to list a  significant Lemma that were originally presented by Brezis and Merle \cite{Brezis} and have also been appeared in \cite{Martinazzi}.

\begin{lem}\label{improve exponential}
	If $\Omega$ is a bounded domain and  $h$ solves
	\begin{align*}
	\begin{cases}
	\displaystyle \left(-\Delta \right)^{m}h=2f\qquad\qquad\qquad\qquad\qquad\qquad\qquad\mathrm{in}\quad  \Omega\\
	\displaystyle h=\left(-\Delta\right)h=\cdots=\left(-\Delta\right)^{m-1}h=0\qquad\qquad\mathrm{on}\quad \partial \Omega,
	\end{cases}
	\end{align*}
	where $f\in L^1(\Omega)$. Then for any $\delta\in (0, nc_n)$, there exist $C_\delta$ such that 
	\begin{align*}
	\int_{\Omega}\exp\left(\frac{\delta |h(x)|}{||f||_{L^1(\Omega)}}\right)\ud x\leq C_{\delta}\left(\mathrm{diam }\, \Omega\right)^n.
	\end{align*}
\end{lem}

\subsection{Asymptotic formula}
In this section, we suppose $(\R^n,e^{2u}|\ud x|^2)$ is complete and  $u$ is a solution of
\begin{align*}
(-\Delta)^{m}u(x)=2Q_g(x)e^{nu(x)}\quad\quad\quad \qquad\mathrm{in}\quad\R^n,
\end{align*}
and  satisfies
\begin{align}\label{sec 5 normal}
		u(x)=\frac{1}{c_n}\int_{\R^n}\log\frac{|y|}{|x-y|}Q_g(y)e^{nu(y)}\ud y+C
\end{align}
and 
\begin{align}\label{sec 5 growth control}
	Q_g(x)\leq C_0|x|^{\gamma}\quad\mathrm{at\quad infinity}\quad \qquad \mathrm{and}\qquad Q_g(x)\geq 0.
\end{align}

IIn the following discussion, we will present the asymptotic formula for solutions \eqref{fine 1} by drawing inspiration from the research conducted by Cheng and Lin \cite{Cheng and Lin} and Martinazzi \cite{Martinazzi}. It is worth noting that similar approaches have been employed in \cite{Hao,Hyder Mancini Martinazzi,Jin and Xiong}.

For any $k\in\N^{+}$, we consider the Green's function of $\left(-\Delta\right)^k$ on $B_4(0)$. Let 
\  \begin{align*}
\begin{cases}
\displaystyle \left(-\Delta\right)^{k}G=\delta_0\qquad\qquad\qquad\qquad\qquad\qquad\qquad\qquad\mathrm{in}\quad  B_{4}(0)\\
\displaystyle G=0,-\Delta G=0,\cdots,\left(-\Delta\right)^{k-1}G=0\qquad\qquad\quad\,\,\,\,\mathrm{on}\quad \partial B_{4}(0),
\end{cases}
\end{align*}
in fact $G$ is a radial function.
For any $\Delta^kh(x)=0$ in $B_8(0)$, then 
\begin{align}\label{formula j}
h(0)=&\int_{B_{4}}\left(-\Delta\right)^{k}Gh=\int_{B_{4}}\left(-\Delta\right)^{k-1}G(-\Delta)h
-\int_{\partial B_{4}}h\frac{\partial}{\partial r}\left(-\Delta\right)^{k-1}G\nonumber\\
=&-\int_{\partial B_{4}}h\frac{\partial}{\partial r}\left(-\Delta\right)^{k-1}G-\int_{B_{4}}(-\Delta)h\frac{\partial}{\partial r}\left(-\Delta\right)^{k-2}G=\cdots\nonumber\\
=&\sum_{i=0}^{k-1}\int_{\partial B_4}\left(-\Delta\right)^ih\left(-\frac{\partial}{\partial r}\right)\left(-\Delta\right)^{k-1-i}G.
\end{align}
Then we claim that $-\frac{\partial}{\partial r}\left(-\Delta\right)^{k-1-i}G=c_i>0$ on $\partial B_4(0)$ for $i=0,1,\cdots,k-1$. Let $h_i$ solves the PDE: 
 \begin{align*}
\begin{cases}
\displaystyle \left(-\Delta\right)^{k}h_i=0\quad\quad\qquad\qquad\qquad\qquad\qquad\qquad\qquad\qquad\mathrm{in}\quad  B_{4}(0)\\
\displaystyle h_i=0,\cdots,\left(-\Delta\right)^i h_i=1,\cdots,\left(-\Delta\right)^{k-1}h_i=0\qquad\quad\,\,\,\,\mathrm{on}\quad \partial B_{4}(0),
\end{cases}
\end{align*}
we know $h_i(x)>0$ in $B_4(0)$ by the maximum principle, then 
\begin{align*}
h_i(0)=\int_{\partial B_{4}}-\frac{\partial}{\partial r}\left(-\Delta\right)^{k-1-i}G=|\partial B_1|c_i>0.
\end{align*}

\begin{lem}\label{bad term}
Suppose $u$ satisfies \eqref{sec 5 normal} and \eqref{sec 5 growth control},
	then for any $\e>0$, there exist $R_{\e}\gg1$ if $|x|>R_{\e}$, 
	\begin{align*}
	\int_{B_1(x)}\log\frac{1}{|x-y|}Q_g(y)e^{nu(y)}\ud y\leq \e\log|x|
	\end{align*}
	and
	\begin{align}\label{fine 1}
	\alpha\log|x|+C\geq v(x)\geq (\alpha-\e)\log|x|-C.
	\end{align}
\end{lem}		
\begin{pf}
	For $|x|>R_0\gg1$, then we solve the PDE
	\begin{align*}
	\begin{cases}
	\displaystyle \left(-\Delta \right)^{m}l=2f\qquad\qquad\qquad\qquad\qquad\qquad\qquad\mathrm{in}\quad  B_4(x)\\
	\displaystyle l=\left(-\Delta\right)l=\cdots=\left(-\Delta\right)^{m-1}l=0\,\qquad\qquad\quad\mathrm{on}\quad \partial B_4(x).
	\end{cases}
	\end{align*}
	By Lemma \ref{improve exponential}, we know there exist $q$ sufficiently big such that
	\begin{align}\label{L1}
	\int_{B_4(x)}e^{nql(y)}\ud y\leq C.
	\end{align}
	Let $h(y)=-(v(y)+l(y))$, then it satisfies 
	\  \begin{align*}
	\begin{cases}
	\displaystyle \left(-\Delta \right)^{m}h=0\qquad\qquad\qquad\qquad\qquad\qquad\qquad\qquad\qquad\qquad\mathrm{in}\quad  B_{4}(x)\\
	\displaystyle h=-v,\left(-\Delta\right)h=\Delta v,\cdots,\left(-\Delta\right)^{m-1}h=-\left(-\Delta\right)^{m-1}v\quad\,\,\,\,\mathrm{on}\quad \partial B_{4}(x).
	\end{cases}
	\end{align*}
\end{pf}
And, let $G$ be the Green's function of $\left(-\Delta\right)^{m-1}$ at $x$, i.e,
\begin{align*}
\begin{cases}
\displaystyle \left(-\Delta\right)^{m-1}G=\delta_x\quad\qquad\qquad\qquad\qquad\qquad\qquad\qquad\qquad\mathrm{in}\quad  B_{4}(x)\\
\displaystyle G=0,-\Delta G=0,\cdots,\left(-\Delta\right)^{m-2}G=0\,\,\,\,\quad\qquad\qquad\qquad\mathrm{on}\quad \partial B_{4}(x).
\end{cases}
\end{align*}
By \eqref{formula j} and Lemma \ref{estimate of sphere integral}, we know 
\begin{align*}
-\Delta h(x)=&\sum_{i=0}^{m-2}\int_{\partial B_4(x)}\left(-\Delta\right)^{i+1}h\left(-\frac{\partial}{\partial r}\right)\left(-\Delta\right)^{m-2-i}G\\
=&\sum_{i=1}^{m-1}c_{i-1}\int_{\partial B_{4}(x)}\left(-\Delta\right)^{i}h=\sum_{i=1}^{m-1}c_{i-1}\int_{\partial B_{4}(x)}-\left(-\Delta\right)^{i}v\leq C.	
\end{align*}
For small $r_0$, $z\in B_{r_0}(x)$, let $G_z$ be the Green's function at $z$ i.e,
\begin{align*}
\begin{cases}
\displaystyle \left(-\Delta\right)^{m-1}G_z=\delta_z\quad\qquad\qquad\qquad\qquad\qquad\qquad\qquad\qquad\mathrm{in}\quad  B_{4}(x)\\
\displaystyle G_z=0,-\Delta G_z=0,\cdots,\left(-\Delta\right)^{m-2}G_z=0\,\,\,\,\qquad\qquad\qquad\mathrm{on}\quad \partial B_{4}(x).
\end{cases}
\end{align*}
Then if $r_0$ is sufficiently small, we know 
\begin{align*}
0<\frac{c_{i}}{2}<-\frac{\partial}{\partial \nu}\left(-\Delta\right)^{m-2-i}G_z<2c_{i}\qquad\mathrm{on}\quad\partial B_{4}(x)
\end{align*} 
for $i=0,1,\cdots,m-2$.
The Green's function at $x$ can be obtained by translating the Green's function at the zero point. Hence, we note that $r_0$ is independent of $x$ ! Similarly, you can get 
\begin{align*}
-\Delta h(z)\leq C\qquad\qquad\qquad\mathrm{in}\quad B_{r_0}(x).
\end{align*}
Using the Lemma \ref{negative part L 1 estimate}, we have
\begin{align*}
\int_{B_4(x)}h^{+}\leq 	\int_{B_4(x)}v^{-}\leq C.
\end{align*}
Standard elliptic estimate implies that
\begin{align*}
\sup_{z\in B_{r_0/2}(x)}h\leq \int_{B_4(x)}h^{+}+C\sup_{z\in B_{r_0}(x)}-\Delta h\leq C.
\end{align*}
By \eqref{sec 5 normal}, we know
\begin{align*}
u=C-v=C+h+l\leq C+l\qquad\mathrm{in}\quad B_{r_0/2}(x).
\end{align*}
Now, for any $q>2$ and $\beta>0$, since $B_{|x|^{-\beta}}(x)\subset B_{r_0/2}(x)$ for $|x|\gg1$,  we conclude that
\begin{align}\label{fine estimate}
\int_{B_{|x|^{-\beta}}(x)}\left(Q_ge^{nu}\right)^{2}\leq& C\int_{B_{|x|^{-\beta}}(x)}Q^2_ge^{2nl}\leq C|x|^{2\gamma }\int_{B_{|x|^{-\beta}}(x)}e^{2nl}\nonumber\\
\leq&C|x|^{2\gamma }\left(\int_{B_{|x|^{-\beta}}(x)}e^{qnl}\right)^{\frac{2}{q}}\left(\int_{B_{|x|^{-\beta}}(x)}\right)^{\frac{q-2}{q}}\nonumber\\
\leq&C|x|^{2\gamma -\frac{\beta n(q-2)}{q}},
\end{align}
where we use \eqref{L1} and choose $\beta=\frac{2\gamma q}{n(q-2)}$. For \eqref{fine estimate}, we obtain
\begin{align}\label{Lem 4.5.1}
\int_{B_1(x)}\log\frac{1}{|x-y|}Q_g(y)e^{nu(y)}\ud y
\leq& \left(\int_{B_{|x|^{-\beta}}(x)}\left(\log\frac{1}{|x-y|}\right)^2\ud y\right)^{\frac{1}{2}}\left(\int_{B_{|x|^{-\beta}}(x)}\left(Q_ge^{nu}\right)^{2}\right)^{\frac{1}{2}}\nonumber\\
+&\int_{B_1(x)\backslash B_{|x|^{-\beta}}(x)}\log\frac{1}{|x-y|}Q_g(y)e^{nu(y)}\ud y\nonumber\\
\leq& C+\beta\log|x|\int_{B_1(x)\backslash B_{|x|^{-\beta}}(x)}Q_g(y)e^{nu(y)}\ud y\nonumber\\
\leq& C+\e\log|x|,
\end{align}
where the last inequality we use the $\int_{\R^n}Q_ge^{nu}<+\infty$ and $|x|\gg1$.	Combining with Lemma \ref{upper bound}, Lemma \ref{lower bound} and \eqref{Lem 4.5.1} , we get the estimate \eqref{fine 1}.

\begin{thm}\label{bonudeness thm}
	Suppose 
	\begin{align*}
	u(x)=\frac{1}{c_n}\int_{\R^n}\log\frac{|y|}{|x-y|}Q_g(y)e^{nu(y)}\ud y+C
	\end{align*}
	and $Q_ge^{nu}$ is absolutely integrable , if 
	\begin{align}
	Q_g(x)\leq C_0|x|^{\gamma}\quad\mathrm{at\quad infinity}\quad \qquad \mathrm{and}\qquad Q_g(x)\geq 0.
	\end{align}
Then, for any $\e>0$, there exist $R_{\e}\gg1$ such that
\begin{align*}
-(\alpha-\e)\log|x|+C\geq u(x)\geq-\alpha\log|x|-C
\end{align*}
for $|x|>R_{\e}$.
	
\end{thm}

\begin{pf}
	
	The asymptotic formula follows by \eqref{upper bounded},  \eqref{lower bounded} and Lemma \ref{bad term}.
\end{pf}

Proof of Theorem \ref{asy thm}: This follows by Theorem \ref{thm A1} and Theorem \ref{bonudeness thm}.\\

\section{Appendix}\label{sec appendix}

Here we give a detail proof above Lemma \ref{upper bound} and Lemma \ref{lower bound}. \\

\textbf{Proof of Lemma \ref{upper bound}}: Recall the assumption, we know $f(y)\geq0$ for $|y|\geq R_0$. We decompose $\R^n$ to $\R^n=A_1\cup A_2\cup A_3$, where 
\begin{align*}
A_1=&\left\{y||y|<|x|/2\right\},\quad A_2=\left\{y||x-y|<|x|/2\right\},\\
A_3=&\left\{y||y|>|x|/2,|x-y|>|x|/2\right\}.
\end{align*}
In $A_2$, we have $|x-y|<\frac{|x|}{2}<|y|$, then $\log\frac{|x-y|}{|y|}\leq0$. Hence, 
\begin{align}\label{lem 2.1 formula 1}
c_nv(x)\leq 	\int_{A_1\cup A_3}\log\frac{|x-y|}{|y|}f(y)\ud y.
\end{align} 
In $A_1$,  we know $\frac{|x|}{2}\leq|x-y|\leq \frac{3|x|}{2}$, hence
\begin{align}\label{lem 2.1 formula 2}
\int_{A_1}\log\frac{|x-y|}{|y|}f(y)\ud y=& \log|x|\int_{A_1}f(y)\ud y+\int_{A_1}\log\frac{|x-y|}{|x|}f(y)\ud y-\int_{A_1}\log|y|f(y)\ud y\nonumber\\
\leq&\log|x|\int_{A_1}f(y)\ud y+C\int_{A_1}|f(y)|\ud y-\int_{|y|<R_0}\log|y| f(y)dy\nonumber\\
\leq&\log|x|\int_{\R^n}f(y)\ud y+C+||f||_{L^{\infty}(B_{R_0})}\int_{|y|<R_0}\left|\log|y|\right|dy.
\end{align}
In $A_3$, $\frac{|y|}{3}\leq |x-y|\leq 3|y|$, so we arrive that
\begin{align}\label{lem 2.1 formula 4}
	\int_{A_3}\log\frac{|x-y|}{|y|}f(y)\ud y\leq C\int_{A_4}f(y)\ud y
\end{align}
So, from \eqref{lem 2.1 formula 1}, \eqref{lem 2.1 formula 2} and \eqref{lem 2.1 formula 4}, we conclude that
\begin{align*}
v(x)\leq\alpha\log|x|+C.
\end{align*}
\\\\

\textbf{Proof of Lemma \ref{lower bound}}:	Similarly we decompose $\R^n$ into $\R^n=A_1\cup A_2\cup A_3\cup A_4$, where \begin{align*}
A_1=&\left\{y||y|<R\right\},\quad\quad A_2=\left\{y||x-y|<|x|/2\right\}\\
A_3=&\left\{y||y|>R,|x-y|>|x|/2,|y|<2|x|\right\},\quad A_4=\left\{y||y|>2|x|\right\},
\end{align*}
where $R>R_0$ to be determined and $|x|\gg R$. Then
\begin{align}\label{a 1}
\int_{\R^n}\log\frac{|x-y|}{\vert y\vert} f(y)\ud y
=&\sum_{i=1}^4 \int_{A_i}\log\frac{|x-y|}{\vert y\vert} f(y)\ud y.\nonumber\\
=&I+II+III+IV.
\end{align}
For $A_1$, we suppose $R\ll |x|$ then 
\begin{align}\label{a 2}
I=&\int_{A_1}\log\frac{|x-y|}{|x|} f(y)\ud y+\log|x|\int_{A_1}f(y)\ud y-\int_{A_1}\log|y|f(y)\ud y\nonumber\\
\geq&-C\int_{A_1} |f(y)|\ud y+\log|x|\int_{A_1}f(y)\ud y-\log R\int_{A_1\cap \{y||y|>1\}}|f(y)|\ud y\nonumber\\
-&||f||_{L^{\infty}(B_{1})}\int_{|y|<1}\log\frac{1}{|y|}dy\nonumber\\
\geq&\log|x|\int_{A_1}f(y)\ud y-C(R).
\end{align}
In $A_2$, $|y|>R>R_0$  and $\frac{1}{2}|x|<|y|<\frac{3}{2}|x|$  we know $f(y)\geq 0$ then 
\begin{align}\label{a 3}
II=&\int_{A_2}\log|x-y|f(y)\ud y-\int_{A_2}\log|y|f(y)\ud y\nonumber\\
\geq&\int_{|y-x|<1}\log|x-y|f(y)\ud y-\left(\int_{A_2}f(y)\ud y\right)\log|x|-C\nonumber\\
\geq&\int_{|y-x|<1}\log|x-y|f(y)\ud y-\left(\int_{|y|>R}f(y)\ud y\right)\log|x|-C.
\end{align}
For $A_3$, then 
\begin{align}\label{a 4}
III=&\int_{A_3}\log|x-y|f(y)\ud y-\int_{A_3}\log|y|f(y)\ud y\nonumber\\
\geq&\left(\int_{A_3}f(y)\ud y\right)\log\left(\frac{|x|}{2}\right)-\left(\int_{A_3}f(y)\ud y \right)\log(2|x|)\nonumber\\
\geq&-C\int_{A_3}f(y)\ud y.
\end{align} 	
For $A_4$, we know $\frac{1}{2}|y|\leq|x-y|\leq\frac{3}{2}|y|$, then
\begin{align}\label{a 5}
IV\geq -C\int_{A_4}f(y)\ud y.
\end{align}	
Combining  \eqref{a 1}, \eqref{a 2}, \eqref{a 3}, \eqref{a 4} with  \eqref{a 5}, we obtain
\begin{align*}
	\frac{1}{c_n}\int_{\R^n}\log\frac{|x-y|}{\vert y\vert} f(y)\ud y
	\geq& \frac{\log|x|}{c_n}\left(\int_{\R^n}f(y)\ud y-2\int_{|y|>R}f(y)\ud y\right)\\
	-&\frac{1}{c_n}\int_{|y-x|<1}\log\frac{1}{|x-y|}f(y)\ud y-C(R).
\end{align*}
For any $\e>0$, by choosing $R=R(\e)$ sufficiently big and letting $|x|\gg R (\e)$, then  we complete the proof.\\\\	

		{\noindent\small{\bf Acknowledgment:} We appreciate the discussions with Postdoctoral Fellow Kai Zhang and Associate Professor Yalong Shi. This research did not receive any specific grant from funding agencies in the public, commercial, or not-for-profit sectors. 
			
	{\noindent\small{\bf Data availability:} Data sharing not applicable to this article as no datasets were generated or analysed during the current study.
	\section*{Declarations}
	{\noindent\small{\bf Conflict of interest:} On behalf of all authors, the corresponding author states that there is no conflict of interest.

		\bibliographystyle{unsrt}

			\bigskip
		
		\noindent S. Zhang
		
		\noindent Department of Mathematics, Nanjing University, \\
		Nanjing 210093, China\\[1mm]
		Email: \textsf{dg21210019@smail.nju.edu.cn}
		
		\medskip  	
		
	\end{document}